\documentclass[10pt]{article}
\usepackage{epsfig,amssymb,amsmath}
\usepackage{palatino}
\usepackage{verbatim}
\usepackage{bm}
\def\PP{{\cal P}}

\newtheorem{lemma}{Lemma}[section]

\newtheorem{theorem}[lemma]{Theorem}

\def\be  {\begin{equation}}
\def\ee  {\end{equation}}
\def\ba  {\begin{eqnarray}}
\def\ea  {\end{eqnarray}}
\def\baa {\begin{eqnarray*}}
\def\eaa {\end{eqnarray*}}
\def\bc  {\begin{comment}}
\def\ec  {\end{comment}}
\makeatletter
\@addtoreset{equation}{section}
\makeatother
\def\proof{\medskip\noindent{\bf Proof.} }
\begin{document}

\title{On the Markov inequality in the $L_2$-norm with the Gegenbauer weight}

\author{A.\,Shadrin, G.\,Nikolov, D.\,Aleksov}

\date{}
\maketitle

%******************************************************************

\begin{abstract}

%******************************************************************
Let $w_{\lambda}(t)=(1-t^2)^{\lambda-1/2}$, $\lambda>-1/2$, be the
Gegenbauer weight function, and $\Vert\cdot\Vert$ denote the
associated $L_2$-norm, i.e.,
$$
\Vert f\Vert:=\Big(\int_{-1}^{1}w_{\lambda}(t)\vert
f(t)\vert^2\,dt\Big)^{1/2}.
$$
Denote by $\PP_n$ the set of algebraic polynomials of degree not
exceeding $n$. We study the best (i.e., the smallest) constant
$c_{n,\lambda}$ in the Markov inequality
$$
\Vert p^{\prime}\Vert\leq c_{n,\lambda}\,\Vert p\Vert,\qquad p\in
\PP_n,
$$
and prove that
$$ c_{n,\lambda}< \frac{(n+1)(n+2\lambda+1)}{2\sqrt{2\lambda+1}},\qquad
\lambda>-1/2\,.
$$
Moreover, we prove that the extremal polynomial in this inequality
is even or odd depending on whether $n$ is even or odd.
\end{abstract}

%******************************************************************

\section{Introduction}
Throughout, $\PP_n$ is the collection of algebraic polynomials of
degree at most $n$. For $\lambda>-1/2$,
$w_{\lambda}(t)=(1-t^2)^{\lambda-1/2}$ is the Gegenbauer weight
function in $[-1,1]$, and $\Vert\cdot\Vert$ is the associated
weighted $L_2[-1,1]$ norm, i.e.,
$$
\Vert
f\Vert:=\Big(\int_{-1}^{1}w_{\lambda}(t)|f(t)|^2\,dt\Big)^{1/2}.
$$
We study the Markov inequality inequality in this norm,

\be\label{e1.1} \Vert f^{\prime}\Vert\leq c\,\Vert f\Vert,\qquad \ \
f\in \PP_n, \ee in particular, we are interested in the sharp
constant $c_{n,\lambda}$ in \eqref{e1.1},

\be\label{e1.2} c_{n,\lambda}:=\sup_{\mathop{}^{f\in\PP_n}_{f\ne
0}}\frac{\Vert f^{\prime}\Vert}{\Vert f\Vert}\,. \ee

A well-known fact (see, e.g., \cite{PD1} or \cite{MMR}) is that the
sharp constant in a $L_2$ Markov-type inequality is equal to the
largest singular value of a certain matrix. Despite of this simple
characterization, not much is known about the sharp constant even in
the classical case of the Gegenbauer weight function. Namely, in the
case of a constant weight ($\lambda=1/2$), E. Schmidt \cite{ES}
found the sharp constant asymptotically, showing that
\begin{equation}\label{e1.3}
c_{n,1/2}=\frac{(n+3/2)^2}{\pi}\,\Big(1-\frac{\pi^2-3}{12(n+3/2)^2}
+\frac{R_n}{(n+3/2)^4}\Big)^{-1}\,, \ \ n> 5
\end{equation}
with $\ -6<R_n<13$.

The situation is better for the $L_2$ norms induced by the other
classical weight functions. For the Hermite weight $w(t)=e^{-t^2},\
t\in \mathbb{R}$ both the sharp constant and the extremal polynomial
are known, resp. $c(n)=\sqrt{2n}$ and $p=H_n$, the $n$-th Hermite
polynomial (see, e.g., \cite{ES}). For the Laguerre weight
$w(t)=t^{\alpha}e^{-t},\ \alpha>-1,\ t\in (0,\infty)$, Tur\'{a}n
\cite{PT} found the sharp constant in the case $\alpha=0$, namely
$$
c(n)=\Big(2\sin\frac{\pi}{4n+2}\Big)^{-1}\,.
$$
For more information on the history of the  $L_2$ Markov-type
inequalities, as well as for more general Markov-type inequalities
for higher order derivatives and with different Laguerre and
Gegenbauer $L_2$ norms, we refer the reader to the recent papers of
A. B\"{o}tcher and P. D\"{o}rfler \cite{BD1, BD2, BD3} and the
references therein. In \cite{BD2} the sharp constants are shown to
be asymptotically equal to the norms of certain Volterra operators.
In general, finding explicitly the best constants in these
Markov-type inequalities and the norms of the related Volterra
operators seem to be equally difficult tasks.

The aim of this paper is to derive an upper bound for the best
constant $c_{n,\lambda}$ in the $L_2$ Markov-type inequality
\eqref{e1.1} associated with the Gegenbauer weight function
$w_{\lambda}$. Our main result is the following:
\begin{theorem}\label{t1.1}
The sharp constant $c_{n,\lambda}$ in the Markov-type inequality in
the $L_2$ norm induced by the Gegenbauer weight function
$w_{\lambda}$, $\ \lambda>-1/2$, satisfies the inequality
$$
c_{n,\lambda}< \frac{(n+1)(n+2\lambda+1)}{2\sqrt{2\lambda+1}} \,.
$$
\end{theorem}

The proof of Theorem \ref{t1.1} is given in the next section.
Section 3 contains some remarks and comments, in particular, we show
that the extremal polynomial in \eqref{e1.2} is even or odd
depending on whether $n$ is even or odd (Theorem 3.1).

%******************************************************************
\section{Proof of Theorem \ref{t1.1}}

Let us start with fixing some notation. For
$\mathbf{t}=(t_1,t_2,\ldots,t_n)\in \mathbb{R}^{n}$, $|\mathbf{t}|$
is the Euclidean norm of $\mathbf{t}$,
$|\mathbf{t}|=(t_1^2+t_2^2+\cdots+t_n^2)^{1/2}$. The unit sphere in
$\mathbb{R}^n$ is denoted by $S^n$, $\ S^n:=\{ \mathbf{t}\in
\mathbb{R}^n\,:\, |t|=1\}$. By $\mathbb{R}^n_{+}$ and $S^n_{+}$ we
denote the subsets of $\mathbb{R}^n$ and $S^n$, respectively, with
non-negative coordinates.\smallskip

The inner product of $\mathbf{x}$ and $\mathbf{y}$,
$\,\mathbf{x},\,\mathbf{y}\in \mathbb{R}^n$, is denoted by
$(\mathbf{x},\mathbf{y})$, thus
$|\mathbf{t}|^2=(\mathbf{t},\mathbf{t})$\,.\smallskip

Finally, the trace of a matrix  $\mathbf{A}=(a_{i,j})_{n\times n}$
is denoted by ${\rm tr}(\mathbf{A})$,
$$
{\rm tr}(\mathbf{A})=\sum_{i=1}^{n}a_{i,i}\,.
$$

The associated with $w_{\lambda}$ orthogonal polynomials are the
Gegenbauer polynomials $\{C_m^{\lambda}(t)\}_{m=0}^{n}$ (called also
ultraspherical polynomials), see. e.g., \cite[Chapther 4.7]{GS}. For
$\lambda\ne 0$ they are normalized so that
$$
\int_{-1}^{1}w_{\lambda}(t)C_j^{\lambda}(t)C_k^{\lambda}(t)\,dt=\delta_{j,k}h_k^2\,,
$$
with $\delta_{j,k}$ being the Kronecker delta and
$$
h_k^2=h_{k,\lambda}^2=\frac{2^{1-2\lambda}\pi
\Gamma(k+2\lambda)}{k!(k+\lambda)\Gamma^2(\lambda)}\,.
$$
Hence, the orthonormal polynomials associated with $w_{\lambda}$ are
given by
$$
p_j(t)=h_j^{-1}\,C_j^{\lambda}(t),\qquad j=0,1,\ldots,n,
$$
where, for simplicity, the dependence on $\lambda$ is omitted
hereafter.

On using the well-known identities (see \cite[eqn. 4.7.11]{GS} and
\cite[eqn. 7.13]{RA})
$$
\frac{d}{dt}\,C_j^{\lambda}(t)=2\lambda\,C_{j-1}^{\lambda+1}(t)\,,
\qquad j\geq 1
$$
$$
C_{j}^{\lambda+1}(t)=\sum_{k=0}^{\lfloor
j/2\rfloor}\frac{j-2k+\lambda}{\lambda}\,C_{j-2k}^{\lambda}(t)\,
$$
(these identities are true for $\lambda\ne 0$ only), we find
\begin{equation}\label{e2.1}
p_j^{\prime}(t)=2\sum_{k=0}^{\lfloor
j-1/2\rfloor}(j-2k-1+\lambda)\,\frac{h_{j-2k-1}}{h_{j}}\,p_{j-2k-1}(t)\,,
\qquad j=1,\ldots,n\,.
\end{equation}
Since the representation of $p_j^{\prime}$ depends on the ratios of
the $h$'s rather than on the $h$'s themselves, we find appropriate
to remove the constant factor in  $h_k$, defining hereafter $h_k$ by
\begin{equation}\label{e2.2}
h_k=h_{k,\lambda}=\Big(\frac{
\Gamma(k+2\lambda)}{k!(k+\lambda)}\Big)^{1/2}\,.
\end{equation}
Notice that now \eqref{e2.1} is true also for $\lambda=0$.\smallskip

Let $p\in\PP_n$ and $\Vert p\Vert=1$. Our goal is to find upper
bound for $\Vert p^{\prime}\Vert$, which, by virtue of \eqref{e1.2},
will be an upper bound for $c_{n,\lambda}$ as well. Without loss of
generality we may assume that
$$
p(t)=\sum_{j=1}^{n}t_j\,p_j(t) \  \mbox{ with }\ \
\mathbf{t}=(t_1,t_2,\ldots,t_n)\in \mathbb{R}^n\,,\ \
|\mathbf{t}|=1\,.
$$

From \eqref{e2.1} we find
\begin{equation}\label{e2.3}
\begin{split}
\Vert p^{\prime}\Vert^2=&4\sum_{k=1}^{\lfloor n/2\rfloor}\Bigg(
\sum_{j=k}^{\lfloor
n/2\rfloor}(2k+\lambda-1)\,\frac{h_{2k-1}}{h_{2j}}\,t_{2j}\Bigg)^2\\
&+4\sum_{k=1}^{\lfloor n+1/2\rfloor}\Bigg(\sum_{j=k}^{\lfloor
n+1/2\rfloor}(2k+\lambda-2)\,\frac{h_{2k-2}}{h_{2j-1}}\,t_{2j-1}\Bigg)^2\\
=&4\big[\big|\mathbf{C}_{\lfloor
n/2\rfloor}\mathbf{t}^{\prime}\big|^2+\big|\mathbf{\widetilde{C}}_{\lfloor
n+1/2\rfloor}\mathbf{t}^{\prime\prime}\big|^2\big]\,.
\end{split}
\end{equation}
Here $\mathbf{t}^{\prime}=(t_2,t_4,\ldots,t_{2\lfloor n/2\rfloor})$,
$\mathbf{t}^{\prime\prime}=(t_1,t_3,\ldots,t_{2\lfloor
n+1/2\rfloor-1})$, and for $m\in \mathbb{N},\ m\geq 2$, the $m\times
m$ matrices $\mathbf{C}_{m}$ and $\mathbf{\widetilde{C}}_{m}$ are
defined by
\begin{equation}\label{e2.4}
\mathbf{C}_m= \begin{pmatrix} \alpha_1\beta_1&\alpha_1\beta_2&\cdot\cdots&\alpha_1\beta_m\\
0&\alpha_2\beta_2&\cdots&\alpha_2\beta_m\\
\vdots&\vdots&\ddots&\vdots\\0&0&\cdots&\alpha_m\beta_m\end{pmatrix}\,,\qquad
\mathbf{\widetilde{C}}_m=
\begin{pmatrix}
\widetilde{\alpha}_1\widetilde{\beta}_1&\widetilde{\alpha}_1\widetilde{\beta}_2
&\cdot\cdots&\widetilde{\alpha}_1\widetilde{\beta}_m\\
0&\widetilde{\alpha}_2\widetilde{\beta}_2&\cdots&\widetilde{\alpha}_2\widetilde{\beta}_m\\
\vdots&\vdots&\ddots&\vdots\\0&0&\cdots&\widetilde{\alpha}_m\widetilde{\beta}_m\end{pmatrix}\,,
\end{equation}
where
\begin{equation}\label{e2.5}
\alpha_k=(2k+\lambda-1)h_{2k-1}=
\Bigg(\frac{\Gamma(2k+2\lambda-1)(2k+\lambda-1)}
{\Gamma(2k)}\Bigg)^{1/2},
\end{equation}
\begin{equation}\label{e2.6}
\beta_k=\frac{1}{h_{2k}}=\Bigg(\frac{\Gamma(2k+1)(2k+\lambda)}{\Gamma(2k+2\lambda)}\Bigg)^{1/2},
\end{equation}
\begin{equation}\label{e2.7}
\widetilde{\alpha}_k=(2k+\lambda-2)h_{2k-2}=\Bigg(\frac{\Gamma(2k+2\lambda-2)(2k+\lambda-2)}
{\Gamma(2k-1)}\Bigg)^{1/2},
\end{equation}
\begin{equation}\label{e2.8}
\widetilde{\beta}_k=\frac{1}{h_{2k-1}}
=\Bigg(\frac{\Gamma(2k)(2k+\lambda-1)}{\Gamma(2k+2\lambda-1)}\Bigg)^{1/2}\,.
\end{equation}
Next, we have
\begin{equation}\label{e2.9}
\begin{split}
\big|\mathbf{C}_{\lfloor
n/2\rfloor}\mathbf{t}^{\prime}\big|^2&=\big(\mathbf{C}_{\lfloor
n/2\rfloor}\mathbf{t}^{\prime},\mathbf{C}_{\lfloor
n/2\rfloor}\mathbf{t}^{\prime}\big)=\big(\mathbf{C}_{\lfloor
n/2\rfloor}^{\top}\mathbf{C}_{\lfloor
n/2\rfloor}\mathbf{t}^{\prime},\mathbf{t}^{\prime}\big)=
\big(\mathbf{A}_{\lfloor
n/2\rfloor}\mathbf{t}^{\prime},\mathbf{t}^{\prime}\big)\\
&\leq\nu_{\lfloor n/2\rfloor}\,|\mathbf{t}^{\prime}|^2\,,
\end{split}
\end{equation}
where $\nu_{\lfloor n/2\rfloor}$ is the largest eigenvalue of the
positive definite matrix
$$
\mathbf{A}_{\lfloor n/2\rfloor}:=\mathbf{C}_{\lfloor
n/2\rfloor}^{\top}\mathbf{C}_{\lfloor n/2\rfloor}\,.
$$

Analogously,
\begin{equation}\label{e2.10}
\big|\mathbf{\widetilde{C}}_{\lfloor
n+1/2\rfloor}\mathbf{t}^{\prime\prime}\big|^2=
\big(\mathbf{\widetilde{A}}_{\lfloor
n+1/2\rfloor}\mathbf{t}^{\prime\prime},\mathbf{t}^{\prime\prime}\big)
\leq\widetilde{\nu}_{\lfloor
n+1/2\rfloor}\,|\mathbf{t}^{\prime\prime}|^2=
\widetilde{\nu}_{\lfloor
n+1/2\rfloor}\,\big(1-|\mathbf{t}^{\prime}|^2\big)\,,
\end{equation}
where $\widetilde{\nu}_{\lfloor n+1/2\rfloor}$ is the largest
eigenvalue of the positive definite matrix
$$
\mathbf{\widetilde{A}}_{\lfloor n+1/2\rfloor}:=
\mathbf{\widetilde{C}}_{\lfloor
n+1/2\rfloor}^{\top}\mathbf{\widetilde{C}}_{\lfloor n+1/2\rfloor}\,.
$$

From \eqref{e2.3}, \eqref{e2.9} and \eqref{e2.10} we deduce that

\begin{equation}\label{e2.11}
\Vert p^{\prime}\Vert^2\leq 4\,\max\big\{\nu_{\lfloor n/2\rfloor},
\widetilde{\nu}_{\lfloor n+1/2\rfloor}\big\}\,.
\end{equation}

From \eqref{e2.4} it follows that
\begin{equation} \label{e2.12}
\mathbf{A}_m=\begin{pmatrix}\alpha_1^2\beta_1^2 &
\alpha_1^2\beta_1\beta_2 & \alpha_1^2\beta_1\beta_3 & \cdots &
\alpha_1^2\beta_1\beta_m\\
\alpha_1^2\beta_1\beta_2 &
\Big(\sum_{j=1}^{2}\alpha_j^2\Big)\beta_2^2 &
\Big(\sum_{j=1}^{2}\alpha_j^2\Big)\beta_2\beta_3 & \cdots &
\Big(\sum_{j=1}^{2}\alpha_j^2\Big)\beta_2\beta_m\\
\alpha_1^2\beta_1\beta_3 &
\Big(\sum_{j=1}^{2}\alpha_j^2\Big)\beta_2\beta_3 &
\Big(\sum_{j=1}^{3}\alpha_j^2\Big)\beta_3^2 & \cdots &
\Big(\sum_{j=1}^{3}\alpha_j^2\Big)\beta_3\beta_m\\ \vdots & \vdots &
\vdots & \ddots & \vdots \\
\alpha_1^2\beta_1\beta_m &
\Big(\sum_{j=1}^{2}\alpha_j^2\Big)\beta_2\beta_m &
\Big(\sum_{j=1}^{3}\alpha_j^3\Big)\beta_3\beta_m & \cdots &
\Big(\sum_{j=1}^{m}\alpha_j^2\Big)\beta_m^2
\end{pmatrix}\,,
\end{equation}
and the same representation holds for $\mathbf{\widetilde{A}}_m$
with $\alpha_k$ and $\beta_k$ replaced by $\widetilde{\alpha}_k$ and
$\widetilde{\beta}_k$, respectively.

Since $\mathbf{A}_m$ and $\mathbf{\widetilde{A}}_m$ are positive
definite matrices, their largest eigenvalues $\nu_m$ and
$\widetilde{\nu}_m$ do not exceed their traces:
\begin{equation}\label{e2.13}
\nu_m\leq {\rm
tr}(\mathbf{A}_m)=\sum_{k=1}^{m}\beta_k^2\sum_{j=1}^{k}\alpha_j^2\,,
\qquad \widetilde{\nu}_m\leq {\rm tr}
(\mathbf{\widetilde{A}}_m)=\sum_{k=1}^{m}\widetilde{\beta}_k^2
\sum_{j=1}^{k}\widetilde{\alpha}_j^2\,.
\end{equation}

For the evaluation of ${\rm tr}(\mathbf{A}_m)$ and ${\rm
tr}(\mathbf{\widetilde{A}}_m)$ we shall need two lemmata.

\begin{lemma}\label{l2.1} The following identities hold:
$$
\begin{array}{lllll}
(i) & \displaystyle{\sum_{j=1}^{k}\alpha_j^2=
\frac{\Gamma(2k+2\lambda+1)}{2(2\lambda+1)\Gamma(2k)}}\,;
  & \qquad &
(ii) & \displaystyle{\beta_k^2\,\sum_{j=1}^{k}\alpha_j^2=
\frac{2k(k+\lambda)(2k+\lambda)}{2\lambda+1}}\,.
\end{array}
$$
\end{lemma}

\proof We apply induction with respect to $k$. For $k=1~$ (i)
reduces to
$$
(\lambda+1)\Gamma(2\lambda+1)=\frac{\Gamma(2\lambda+3)}{2(2\lambda+1)}\,,
$$
which easily follows from the property  $\Gamma(x+1)=x\,\Gamma(x)$.
Assume that (i) is true for $k-1\in \mathbb{N}$, then  the induction
hypothesis and \eqref{e2.5} imply
\[
\begin{split}
\sum_{j=1}^{k}\alpha_j^2&=(2k+\lambda-1)\,
\frac{\Gamma(2k+2\lambda-1)}{\Gamma(2k)}
+\frac{\Gamma(2k+2\lambda-1)}{2(2\lambda+1)\Gamma(2k-2)}
\\
&=\frac{\Gamma(2k+2\lambda-1)}{2(2\lambda+1)\Gamma(2k)}\,
\big[2(2\lambda+1)(2k+\lambda-1)+(2k-2)(2k-1)\big]\,.
\end{split}
\]
For the term in the square brackets we have
$$
2(2\lambda+1)(2k+\lambda-1)+(2k-2)(2k-1)=(2k+2\lambda-1)(2k+2\lambda)\,,
$$
whence
$$
\sum_{j=1}^{k}\alpha_j^2=(2k+2\lambda-1)(2k+2\lambda)\,
\frac{\Gamma(2k+2\lambda-1)}{2(2\lambda+1)\Gamma(2k)}
=\frac{\Gamma(2k+2\lambda+1)}{2(2\lambda+1)\Gamma(2k)}\,.
$$
The induction step is done, and this proves claim (i).

Now (ii) follows from (i) and \eqref{e2.6}:
$$
\beta_k^2\,\sum_{j=1}^{k}\alpha_j^2=
\frac{\Gamma(2k+2\lambda+1)}{2(2\lambda+1)\Gamma(2k)}\cdot
\frac{(2k+\lambda)\Gamma(2k+1)}{\Gamma(2k+2\lambda)}=
\frac{2k(k+\lambda)(2k+\lambda)}{2\lambda+1}\,.%\eqno($\Box$)
$$
\begin{lemma}\label{l2.2} The following identities hold:
$$
\begin{array}{lllll}
(i) & \displaystyle{\sum_{j=1}^{k}\widetilde{\alpha}_j^2=
\frac{\Gamma(2k+2\lambda)}{2(2\lambda+1)\Gamma(2k-1)}}\,;
  & \quad &
(ii) &
\displaystyle{\widetilde{\beta}_k^2\,\sum_{j=1}^{k}\widetilde{\alpha}_j^2=
\frac{(2k-1)(2k+\lambda-1)(2k+2\lambda-1)}{2(2\lambda+1)}}\,.
\end{array}
$$
\end{lemma}

\proof For the proof of (i) we use induction with respect to $k$.
For $k=1~$  claim (i) becomes
$$
\lambda\Gamma(2\lambda)=\frac{\Gamma(2\lambda+2)}{2(2\lambda+1)}\,,
$$
and it is obviously true. If we assume that (i) is true for $k-1\in
\mathbb{N}$, then our assumption and \eqref{e2.7} imply
\[
\begin{split}
\sum_{j=1}^{k}\widetilde{\alpha}_j^2&=
(2k+\lambda-2)\,\frac{\Gamma(2k+2\lambda-2)}{\Gamma(2k-1)}
+\frac{\Gamma(2k+2\lambda-2)}{2(2\lambda+1)\Gamma(2k-3)}
\\
&=\frac{\Gamma(2k+2\lambda-2)}{2(2\lambda+1)\Gamma(2k-1)}\,
\big[2(2\lambda+1)(2k+\lambda-2)+(2k-3)(2k-2)\big]\,.
\end{split}
\]
It is easily seen that
$$
2(2\lambda+1)(2k+\lambda-2)+(2k-3)(2k-2)=(2k+2\lambda-2)(2k+2\lambda-1)\,.
$$
Hence,
$$
\sum_{j=1}^{k}\widetilde{\alpha}_j^2=(2k+2\lambda-2)(2k+2\lambda-1)\,
\frac{\Gamma(2k+2\lambda-2)}{2(2\lambda+1)\Gamma(2k-1)}=
\frac{\Gamma(2k+2\lambda)}{2(2\lambda+1)\Gamma(2k-1)}\,,
$$
and the induction proof of (i) is accomplished. The proof of (ii)
follows from (i) and \eqref{e2.8}:
$$
\widetilde{\beta}_k^2\,\sum_{j=1}^{k}\widetilde{\alpha}_j^2=
\frac{\Gamma(2k+2\lambda)}{2(2\lambda+1)\Gamma(2k-1)}\cdot
\frac{(2k+\lambda-1)\Gamma(2k)}{\Gamma(2k+2\lambda-1)}=
\frac{(2k-1)(2k+\lambda-1)(2k+2\lambda-1)}{2(2\lambda+1)}\,.
$$

We are prepared to evaluate the traces of $\mathbf{A}_m$ and
$\mathbf{\widetilde{A}}_m$.

\begin{lemma}\label{l2.3}
The traces of the matrices $\mathbf{A}_m$ and
$\mathbf{\widetilde{A}}_m$ are given by
\begin{equation}\label{e2.14}
{\rm tr}(\mathbf{A}_m)
=\frac{m(m+1)(m+\lambda)(m+\lambda+1)}{2\lambda+1}\,,
\end{equation}
\begin{equation}\label{e2.15}
{\rm tr}(\mathbf{\widetilde{A}}_m)
=\frac{m(m+\lambda)(m^2+\lambda\,m-1/2)}{2\lambda+1}\,.
\end{equation}
Moreover, for every $m\in \mathbb{N}$,
\begin{equation}\label{e2.16}
{\rm tr}(\mathbf{\widetilde{A}}_m)<{\rm tr}(\mathbf{A}_m)< {\rm
tr}(\mathbf{\widetilde{A}}_{m+1})\,.
\end{equation}
\end{lemma}

\proof We start with the proof of \eqref{e2.14}. By \eqref{e2.13}
and Lemma \ref{l2.1} (ii) we obtain
$$
{\rm
tr}(\mathbf{A}_m)=\sum_{k=1}^{m}\beta_k^2\Big(\sum_{j=1}^{k}\alpha_j^2\Big)
=\frac{2}{2\lambda+1}\,\sum_{k=1}^{m}k(k+\lambda)(2k+\lambda)\,.
$$
The proof of \eqref{e2.14} is accomplished by showing by induction
that
$$
2\sum_{k=1}^{m}k(k+\lambda)(2k+\lambda)
=m(m+1)(m+\lambda)(m+\lambda+1)\,.
$$
Indeed, the latter identity is true for $m=1$, and assuming it is
true for $m-1\in \mathbb{N}$, we get \baa
&&2\sum_{k=1}^{m}k(k+\lambda)(2k+\lambda)
=(m-1)m(m+\lambda-1)(m+\lambda)+ 2m(m+\lambda)(2m+\lambda)\\
&&=m(m+\lambda)\big[(m-1)(m+\lambda-1)+2(2m+\lambda)\big] =
m(m+1)(m+\lambda)(m+\lambda+1)\,.\eaa

The proof of \eqref{e2.15} is similar to that of \eqref{e2.14}. We
make use of \eqref{e2.13} and Lemma \ref{l2.2} (ii) to obtain
$$
{\rm
tr}(\mathbf{\widetilde{A}}_m)=\sum_{k=1}^{m}\widetilde{\beta}_k^2
\Big(\sum_{j=1}^{k}\widetilde{\alpha}_j^2\Big)
=\frac{1}{2(2\lambda+1)}\,\sum_{k=1}^{m}(2k-1)(2k+\lambda-1)(2k+2\lambda-1)\,,
$$
and then apply induction with respect to $m$ to show that
$$
\sum_{k=1}^{m}(2k-1)(2k+\lambda-1)(2k+2\lambda-1)=m(m+\lambda)(2m^2+2\lambda\,m-1)\,.
$$
This identity is obviously true for $m=1$, and if it is true for
$m-1\in \mathbb{N}$, then an easy calculation yields
\[
\begin{split}
\sum_{k=1}^{m}(2k-1)(2k+\lambda-1)(2k+2\lambda-1)
=&(m-1)(m+\lambda-1)\big(2(m-1)^2+2\lambda(m-1)-1\big))\\
&+(2m-1)(2m+\lambda-1)(2m+2\lambda-1))\\
=&m(m+\lambda)(2m^2+2\lambda m-1)\,.
\end{split}
\]
The proof of \eqref{e2.15} is accomplished. The verification of
inequalities \eqref{e2.16} is straightforward, and therefore is
omitted. Lemma \ref{l2.3} is proved.
\medskip

If $n$ is even, $n=2m$, then $\lfloor n/2\rfloor = \lfloor
n+1/2\rfloor =m$, and by \eqref{e2.11}, \eqref{e2.13} and Lemma
\ref{l2.3} we obtain
\[
\begin{split}
\Vert p^{\prime}\Vert^2\leq &4\max\{\nu_m,\widetilde{\nu}_m\}\leq
4\max\{{\rm tr}(\mathbf{A}_m),{\rm
tr}(\mathbf{\widetilde{A}}_m)\}=4\,{\rm
tr}(\mathbf{A}_m)\\
=&\frac{4m(m+1)(m+\lambda)(m+\lambda+1)}{2\lambda+1}=
\frac{n(n+2)(n+2\lambda)(n+2\lambda+2)}{4(2\lambda+1)}\\
<&\frac{(n+1)^2(n+2\lambda+1)^2}{4(2\lambda+1)}\,.
\end{split}
\]

If $n$ is odd, $n=2m+1$, then $\lfloor n/2\rfloor = m$ and $\lfloor
n+1/2\rfloor =m+1$. From \eqref{e2.11}, \eqref{e2.13} and
Lemma~\ref{l2.3} we deduce
\[
\begin{split}
\Vert p^{\prime}\Vert^2\leq&4\max\{\nu_m,\widetilde{\nu}_{m+1}\}\leq
4\max\{{\rm tr}(\mathbf{A}_m),{\rm
tr}(\mathbf{\widetilde{A}}_{m+1})\}=4\,{\rm
tr}(\mathbf{\widetilde{A}}_{m+1})\\
=&\frac{4(m+1)(m+\lambda+1)\big((m+1)^2+\lambda(m+1)-1/2\big)}{2\lambda+1}\\
=&\frac{(n+1)(n+2\lambda+1)\big((n+1)^2+2\lambda(n+1)-2\big)}{4(2\lambda+1)}\\
<&\frac{(n+1)^2(n+2\lambda+1)^2}{4(2\lambda+1)}\,.
\end{split}
\]
Thus, whenever $p\in\PP_n$ and $\Vert p\Vert=1$, we have
$$
\Vert p^{\prime}\Vert<
\frac{(n+1)(n+2\lambda+1)}{2\sqrt{2\lambda+1}}\,
$$
or, equivalently,
\begin{equation}\label{e2.17}
\Vert p^{\prime}\Vert<
\frac{(n+1)(n+2\lambda+1)}{2\sqrt{2\lambda+1}}\,\Vert
p\Vert\,,\qquad p\in\PP_n\,.
\end{equation}
The latter inequality proves the desired estimate for
$c_{n,\lambda}$. The proof of Theorem \ref{t1.1} is complete.
%******************************************************************
\section{Remarks}
\textbf{(1)}~~ Our proof of Theorem \ref{t1.1} makes use of the fact
that the traces of the positive definite matrices
$\mathbf{\widetilde{A}}_{m}$, $\mathbf{A}_m$, and
$\mathbf{\widetilde{A}}_{m+1}$ satisfy inequalities \eqref{e2.16}.
We show below that, in fact, the largest eigenvalues of these
matrices are ordered in the same way.
\begin{lemma}\label{l3.1}
The largest eigenvalues $\widetilde{\nu}_{m},\ \nu_{m}$ and
$\widetilde{\nu}_{m+1}$ of the matrices
$\mathbf{\widetilde{A}}_{m}$, $\mathbf{A}_m$ and
$\mathbf{\widetilde{A}}_{m+1}$, respectively, satisfy the
inequalities
\begin{equation}\label{e3.1}
(0<)\ \ \widetilde{\nu}_{m}<\nu_{m}<\widetilde{\nu}_{m+1}\,.
\end{equation}
\end{lemma}

\proof Obviously, the elements $a_{k,i}$ of matrices $\mathbf{A}_m$
and the elements $\widetilde{a}_{k,i}$ of matrices
$\mathbf{\widetilde{A}}_m$ are positive, see see \eqref{e2.12} and
\eqref{e2.5} -- \eqref{e2.8}. We shall prove that they satisfy the
inequalities
\begin{equation}\label{e3.2}
\widetilde{a}_{k,i}<a_{k,i}<\widetilde{a}_{k+1,i+1}\,.
\end{equation}
Since $\mathbf{A}_m$ and $\mathbf{\widetilde{A}}_m$ are symmetric
matrices, we may assume that $i\geq k\geq 1$\,, in which case
(see\eqref{e2.12})
$$
a_{k,i}=\frac{\beta_i}{\beta_k}\,\Big(\beta_k^2\,\sum_{j=1}^{k}\alpha_j^2\Big)\,,
\qquad
\widetilde{a}_{k,i}=\frac{\widetilde{\beta}_i}{\widetilde{\beta}_k}\,
\Big(\widetilde{\beta}_k^2\,\sum_{j=1}^{k}\widetilde{\alpha}_j^2\Big)\,.
$$
Let us prove first the left inequality in \eqref{e3.2}. On using
\eqref{e2.6}, \eqref{e2.8} and Lemmas \ref{l2.1}, \ref{l2.2}, we
obtain
\begin{equation}\label{e3.3}
\frac{a_{k,i}}{\widetilde{a}_{k,i}}=\frac{2k+2\lambda}{2k-1}\,
\sqrt{\varphi(2k)}\,\sqrt{\varphi(2i)}\,,
\end{equation}
where
$$
\varphi(x):=\frac{x(x+\lambda)}{(x+\lambda-1)(x+2\lambda-1)}\,,\qquad
x\geq 2\,.
$$
We consider separately two cases:\smallskip

{\it Case (i): $-1/2<\lambda\leq 1$}.~~ In this case $\varphi(x)\geq
1$ for $x\geq 2$, and since $(2k+2\lambda)/(2k-1)>1$, we conclude
from \eqref{e3.3} that $a_{k,i}/\widetilde{a}_{k,i}>1$.\smallskip

{\it Case (ii): $\lambda> 1$}.~~Since
$$
\varphi^{\prime}(x)=\frac{\lambda-1}{\lambda}\,
\Big(\frac{2\lambda-1}{(x+2\lambda-1)^2}+\frac{1}{(x+\lambda-1)^2}\Big)
\,,
$$
in this case $\varphi(x)$ is monotone increasing for $x\geq 2$. From
\eqref{e3.3} and  $i\geq k$ we obtain
$$
\frac{a_{k,i}}{\widetilde{a}_{k,i}}\geq
\frac{2k+2\lambda}{2k-1}\,\varphi(2k)=\frac{2k(2k+\lambda)(2k+2\lambda)}
{(2k-1)(2k+\lambda-1)(2k+2\lambda-1)}>1\,.
$$

The right inequality in \eqref{e3.2} is proved in the same fashion.
Using again \eqref{e2.6}, \eqref{e2.8} and Lemmas \ref{l2.1},
\ref{l2.2}, we obtain with the same $\varphi$ as above that
$$
\frac{\widetilde{a}_{k+1,i+1}}{a_{k,i}}=\frac{2k+2\lambda+1}{2k}\,
\sqrt{\varphi(2k+1)}\,\sqrt{\varphi(2i+1)}\,.
$$
The right-hand side here coincides with the right-hand side of
\eqref{e3.3} with $k$ and $i$ replaced by $k+1/2$ and $i+1/2$,
respectively, hence it is greater than 1. Inequalities \eqref{e3.2}
are proved.\smallskip

We are ready now to prove inequalities \eqref{e3.1}. Let
$\widetilde{\mathbf{t}}\in S^{m}_{+}$ be the eigenvector of
$\mathbf{\widetilde{A}}_m$ corresponding to the eigenvalue
$\widetilde{\nu}_{m}$, then by \eqref{e3.2}
\begin{equation}\label{e3.4}
\begin{split}
\widetilde{\nu}_{m}&=\sup_{\mathbf{t}\in S^{m}}
(\mathbf{\widetilde{A}}_m\mathbf{t},\mathbf{t}) =\sup_{\mathbf{t}\in
S^{m}_{+}}\,(\mathbf{\widetilde{A}}_m\mathbf{t},\mathbf{t})=
(\mathbf{\widetilde{A}}_m\widetilde{\mathbf{t}},\widetilde{\mathbf{t}})\\
&< (\mathbf{A}_m\widetilde{\mathbf{t}},\widetilde{\mathbf{t}})\leq
\sup_{\mathbf{t}\in S^{m}_{+}}(\mathbf{A}_m\mathbf{t},\mathbf{t})=
\sup_{\mathbf{t}\in
S^{m}}(\mathbf{A}_m\mathbf{t},\mathbf{t})={\nu}_{m}\,,
\end{split}
\end{equation}
whence the first inequality in \eqref{e3.1} is established.
\smallskip

If $\bm{\tau}\in S^{m}_{+}$ is the eigenvector of $\mathbf{A}_m$
corresponding to the eigenvalue $\nu_{m}$\,,
$\bm{\widetilde{\tau}}=(0,\bm{\tau})\in S^{m+1}_{+}$, and
$\mathbf{\hat{A}}_m$ is the $m\times m$ matrix obtained from
$\mathbf{\widetilde{A}}_{m+1}$ by deletion of its first row and
column, then by \eqref{e3.2}
\begin{equation}\label{e3.5}
\begin{split}
\nu_m&=\sup_{\mathbf{t}\in
S^{m}}(\mathbf{A}_m\mathbf{t},\mathbf{t})=\sup_{\mathbf{t}\in
S^{m}_{+}}(\mathbf{A}_m\mathbf{t},\mathbf{t})=(\mathbf{A}_m\bm{\tau},
\bm{\tau})\\
&<(\mathbf{\hat{A}}_m\bm{\tau}, \bm{\tau})
=(\mathbf{\widetilde{A}}_{m+1}\widetilde{\bm{\tau}},\widetilde{\bm{\tau}})
\leq \sup_{\mathbf{t}\in
S^{m+1}_{+}}(\mathbf{\widetilde{A}}_{m+1}\mathbf{t},\mathbf{t}) =
\sup_{\mathbf{t}\in
S^{m+1}}(\mathbf{\widetilde{A}}_{m+1}\mathbf{t},\mathbf{t})=\widetilde{\nu}_{m+1}\,,
\end{split}
\end{equation}
and the second inequality in \eqref{e3.1} is settled. The proof of
Lemma 3.1 is complete.\medskip

\textbf{(2)}~~With Lemma \ref{l3.1} at our disposal, let us turn
back to the proof of Theorem \ref{t1.1}. If $p\in\PP_n$ and $\Vert
p\Vert=1$, it follows from \eqref{e2.3}, \eqref{e2.9}, \eqref{e2.10}
combined with \eqref{e3.4} and \eqref{e3.5} that the upper bound for
$\Vert p^{\prime}\Vert^2$ given by \eqref{e2.11} is attainable, and
this upper bound $ $ is the best possible. If, e.g., $n=2m$, then,
with the notation from the proofs of Theorem \ref{t1.1} and Lemma
\ref{l3.1} we have
$$
\Vert p^{\prime}\Vert^2=4\nu_m
$$
provided $\mathbf{t}^{\prime}=\bm{\tau}\in S^m_{+}$ is the
eigenvector of $\mathbf{A}_m$ corresponding to its largest
eigenvalue $\nu_m$ and $\mathbf{t}^{\prime\prime}$ is the null
vector in $\mathbb{R}^{m}$. Analogously, if $n=2m+1$, then
$$
\Vert p^{\prime}\Vert^2=4\widetilde{\nu}_{m+1}
$$
provided $\mathbf{t}^{\prime}$ is the null vector in
$\mathbb{R}^{m}$ and $\mathbf{t}^{\prime\prime}\in S^{m+1}_{+}$ is
the eigenvector of $\widetilde{\mathbf{A}}_{m+1}$ corresponding to
$\widetilde{\nu}_{m+1}$, the largest eigenvalue of
$\widetilde{\mathbf{A}}_{m+1}$.\smallskip

Thus, we obtain the following (rather natural) result, which, by a
different argument, has been proved for $\lambda\geq 0$ in
\cite{GN}:
\begin{theorem}\label{t3.2}
The sharp constant $c_{n,\lambda}$ in the Markov inequality in the
$L_2$-norm induced by the Gegenbauer weight $w_{\lambda}$,
$\lambda>-1/2$, is given by
$$
c_{n,\lambda}=\begin{cases}2\sqrt{\nu_m}\,,&\ \mbox{ if }\ n=2m,\\
2\sqrt{\widetilde{\nu}_{m+1}}\,,&\ \mbox{ if }\ n=2m+1\,.\end{cases}
$$
Moreover, if $p\in\PP_n$, $p\ne 0$,  is an extremal polynomial in
this inequality, i.e., $\Vert p^{\prime}\Vert=c_{n,\lambda}\,\Vert
p\Vert$\,, then, with some $\theta\in \mathbb{R}$,
$$
p(t)=e^{i\theta}\sum_{k=0}^{\lfloor
n-1/2\rfloor}t_{n-2k}\,C_{n-2k}^{\lambda}(t)\,,\quad \mbox{ with }\
\ t_n>0\quad \mbox{ and }\ t_{n-2k}\geq 0,\ \ 1\leq k\leq
\Big\lfloor \frac{n-1}{2}\Big\rfloor\,.
$$
In particular, $p$ is an even (resp., odd) polynomial if $n$ is even
(resp., odd).
\end{theorem}
\medskip

\textbf{(3)}~~ A brief look at the final part of the proof of
Theorem \ref{t1.1} shows that we can get slightly better upper
bounds for $\Vert p^{\prime}\Vert$, and hence for $c_{n,\lambda}$,
if we distinguish between the cases of even and odd $n$. As the
resulting improvement is neglectable, and also requires separation
of the cases $-1/2<\lambda<1/2~$ and $~\lambda\geq 1/2$, we skip the
details.
\medskip

\textbf{(4)}~~ As a special case of \cite[Theorem 1.2]{BD3} (which
can be obtained also by a result of D\"{o}rfler \cite{PD2} about the
best constants in the $L_2$ Markov inequalities with Laguerre's
weights combined with an asymptotic relation between the best
constants in the $L_2$ Markov inequalities with the Laguerre and the
Gegenbauer weights \cite[Theorem 1.2]{BD1}), we have
\begin{equation}\label{e3.6}
\lim_{n\rightarrow\infty}\frac{c_{n,\lambda}}{n^2}
=\frac{1}{2\,j_{\frac{2\lambda-3}{4},1}}\,,
\end{equation}
where $j_{\nu,1}$ is the first positive zero of the Bessel function
$J_{\nu}$. Hence, $c_{n,\lambda}=O(n^2)$ as $n\rightarrow\infty$,
and our upper bound for $c_{n,\lambda}$ given by Theorem \ref{t1.1}
has the right order with respect to $n$. Let us point out that our
Markov inequality \eqref{e2.17} holds for all $\lambda>-1/2$ and
$n\in \mathbb{N}$. In contrast, \eqref{e3.6} implies that whenever
$c>\frac{1}{2\,j_{\frac{2\lambda-3}{4},1}}$\,, we have $\Vert
p^{\prime}\Vert\leq c\,n^2\,\Vert p\Vert$ for every $p\in \PP_n$, if
$n\geq n_0(c,\lambda)$. However, in general, neither $j_{\nu,1}$ nor
$n_0(c,\lambda)$ is known explicitly.
\smallskip

By another result of D\"{o}rfler \cite{PD1}, (see \cite[eqn.
(4)]{PD2}), we have
\begin{equation}\label{e3.7}
\frac{1}{\sqrt{2(2\lambda+1)(2\lambda+5)}}
\leq\lim_{n\rightarrow\infty}\frac{c_{n,\lambda}}{n^2}\leq
\frac{1}{2\sqrt{2\lambda+1}}\,,
\end{equation}
thus Theorem \ref{t1.1} furnishes the same upper bound for
$\lim\,c_{n,\lambda}/n^2$ as \eqref{e3.7}. A comparison of Schmidt's
result \cite{ES} for $\lambda=1/2$ and the upper bound in
\eqref{e3.7} shows that
$$
\lim_{n\rightarrow\infty}\frac{c_{n,1/2}}{n^2}
=\frac{1}{\pi}=0.318309\ldots<\frac{1}{2\sqrt{2}}=0.353553\,.
$$
In the special cases $\lambda=0$ and $\lambda=1$ Nikolov \cite{GN}
proved the following tighter estimates:
$$
0.472135 \leq\lim_{n\rightarrow\infty}\frac{c_{n,0}}{n^2}\leq
0.478849\,,\ \ 0.248549
\leq\lim_{n\rightarrow\infty}\frac{c_{n,1}}{n^2}\leq 0.256861\,.
$$

\medskip
{\bf Acknowledgements.} The second- and the third named authors are
partially supported by the Sofia University Research Fund through
Grant 75/2015.

\bigskip\bigskip\noindent
%======================================================================
{\sc Alexei Shadrin} \smallskip\\
Department of Applied Mathematics and Theoretical Physics (DAMTP) \\
Cambridge University \\
Wilberforce Road \\
Cambridge CB3 0WA \\
UNITED KINGDOM \\
{\it E-mail:} {\tt a.shadrin@damtp.cam.ac.uk}
\bigskip

\bigskip\noindent
{\sc Geno Nikolov, Dragomir Aleksov} \smallskip\\
Department of Mathematics \\
Universlty of Sofia \\
5 James Bourchier Blvd. \\
1164 Sofia \\
BULGARIA \\
{\it E-mails:}\hspace*{-6mm}
\begin{tabular}{ll}
&{\tt geno@fmi.uni-sofia.bg}\vspace*{-1ex}\\
&{\tt dragomira@fmi.uni-sofia.bg}
\end{tabular}
\end{document}